\documentclass[12]{article}
\usepackage[dvips]{graphicx}
\def\no{\noindent}
\def\beq{\begin{equation}}
\def\eeq{\end{equation}}

\def\h{{\ \buildrel \rm h \over = \ }}
\def\R{{I\!\!R}}

\def\no{\noindent}
\def\beq{\begin{equation}}
\def\eeq{\end{equation}}

\begin{document}
\title{Generalized Taylor's Theorem}

\author{Garret Sobczyk, \\
Departamento de Actuar\'ia, F\'isica y Matem\'aticas \\ Universidad de Las Am\'ericas - Puebla,\\ 72820 Cholula, Mexico}
\maketitle

\begin{abstract} The Euclidean algorithm makes possible a simple but powerful generalization of
Taylor's theorem. Instead of expanding a function in a series around a single point, one
spreads out the spectrum to include any number of points with given multiplicities. Taken together
with a simple expression for the remainder, this theorem becomes a powerful
tool for approximation and interpolation in numerical analysis. 
We also have a corresponding theorem for rational approximation.  
   \end{abstract}
\section*{Generalize Rolle's Theorem}
   Let $h(x)=\prod_{i=1}^r (x-x_i)^{m_i}$ for distinct $x_i\in [a,b]\subset \R$
    with multiplicity $m_i \ge 1$, 
   and let $n=deg(h(x))$. Given two functions $f(x)$ and $g(x)$, we say that
    $f(x)=g(x)\ mod(h(x))$ or $f(x) \h  g(x)$ if for each $1\le i \le r$ and $0 \le k < m_i$
   \beq  f^{(k)}(x_i)=g^{(k)}(x_i). \label{conditions} \eeq
If $f(x)$ and $g(x)$ are polynomials, then $f(x)\h g(x)$ is equivalent to saying
that if $f(x)$ and $g(x)$ are divided by the polynomial $h(x)$ (the Euclidean algorithm), they
give the same remainder.
 We denote the {\it factor ring of polynomials} modulo $h(x)$ over the real numbers $\R$ by
\[ \R[x]_h :=\R[x]/<h(x)>, \]
see \cite[p.266]{Gal}.

\bigskip

{\bf Generalized Rolle's Theorem: }
   Let $f(x)\in C[a,b]$ and $(n-1)$-times differentiable on $(a,b)$.
   If $f(x)=0 \ mod(h(x))$ , then there exist a $c\in (a,b)$ such that $f^{(n-1)}(c)=0$.
      
   {\bf Proof:} Following \cite[p.38]{Linz},
define the function $\sigma(u,v):=\Big(\begin{array}{cc} 1, & u<v \cr 0, & u\ge v \end{array}\Big)$.
   The function $\sigma$ is needed to count the {\it simple zeros} of the polynomial $h(x)$ and its derivatives.

     Let
   $\# h^{(k)}$ denote the number of simple zeros that the polynomial equation 
$h^{(k)}(x)=0$ has. Clearly $\# h = r = \sum_{i=1}^r\sigma(0,m_i).$ By the Classical Rolle's theorem
  \[ \# h^\prime = \sum_{i=1}^r \sigma(1,m_i)+ (\# h)-1= \sum_{i=1}^r \sigma(1,m_i)+\sum_{i=1}^r\sigma(0,m_i)-1. \] 
Continuing this process, we find that 
 \[ \# h^{\prime \prime}=\sum_{i=1}^r \sigma(2,m_i)+(\# h^\prime)-1=\sum_{i=1}^r \sigma(2,m_i)+ 
\sum_{i=1}^r \sigma(1,m_i)+\sum_{i=1}^r\sigma(0,m_i)-2,\]
 and more generally that  
  \[ \# h^{(k)}=\sum_{i=1}^r \sigma(k,m_i)+(\# h^{(k-1)})-1 = 
         \sum_{i=1}^r \sum_{j=0}^k \sigma(j,m_i) - k  \]
for all integers $k\ge 0$.
  
  For $k=n-1$, we have
  \[ \# h^{(n-1)}= \sum_{i=1}^r \sum_{j=0}^{n-1} \sigma(j,m_i) - (n-1) =\sum_{i=1}^r m_i - (n-1)=1 . \]
   The proof is completed by noting that $\# f^{(k)}\ge \# h^{(k)}$ for each $k\ge 0$, and hence
$\# f^{(n-1)}\ge 1$. 
   
    \hfill $ \framebox{} $        
       
 \section*{Approximation Theorems:} 
 
 We can now prove 
 
 {\bf Generalized Taylor's Theorem:}
   Let $f(x)\in C[a,b]$ and $n$ times differentiable on $(a,b)$. 
   Suppose that $f(x)= g(x)\ mod(h(x))$ for some polynomial $g(x)$ where $deg(g)<deg(h)$. 
Then for every $x\in [a,b]$ there exist a $c\in (a,b)$ such that
     \[ f(x)=g(x)+\frac{f^{(n)}(c)}{n!}h(x). \]
     
        {\bf Proof:} For a given $x\in [a,b]$, define the function
       \beq p(t)=f(t)-g(t)-[f(x)-g(x)]\frac{h(t)}{h(x)}. \label{*} \eeq
    In the case that $x=x_i$ for some $1 \le i \le r$, it is shown below that $p(t)$ has a
removable singularity and can be redefined accordingly. Noting that
      \[  p(t)=0 \ mod(h(t)) \quad {\rm and} \quad p(x)=0 \]
    it follows that $p(t) = 0 \ mod(h(t)(t-x))$. Applying
    the Generalized Rolle's Theorem to $p(t)$, there exists
    a $c\in (a,b)$ such that $p^{(n)}(c)=0$. Using (\ref{*}), we calculate  $p^{(n)}(t)$, getting
      \[   p^{(n)}(t)=f^{(n)}(t)-g^{(n)}(t)-[f(x)-g(x)]\Big(\frac{d}{dt}\Big)^n \frac{h(t)}{h(x)} \]
     \[ = f^{(n)}(t)-[f(x)-g(x)]\frac{n!}{h(x)}, \]
   so that
       \[ 0=p^{(n)}(c)= f^{(n)}(c)-[f(x)-g(x)]\frac{n!}{h(x)} \] 
   from which the result follows.
   
    \hfill $ \framebox{} $        
 
\no    Applying the theorem to the case when $x=x_i$, we find by repeated application of L'Hospital's rule that
   \[ \lim_{x\rightarrow x_i} \frac{f(x)-g(x)}{h(x)}=\frac{f^{(m_i)}(x_i)-g^{(m_i)}(x_i)}{h^{(m_i)}(x_i)}
                              =\frac{f^{(n)}(c)}{n!}. \] 

There remains the question of how do we calculate the polynomial $g(x)$ with the property that
$f(x)\h g(x)$ where $deg(g(x))< deg(h(x))$? The {\it brute force} method is to impose the conditions
(\ref{conditions}), and solve the resulting system of linear equations for the unique solution known as the
{\it osculating polynomial} approximation to $f(x)$, see \cite{Davis}, \cite[p.52]{Sto}. A far more
powerful method is to make use of the special algebraic properties of the {\it spectral basis} of the factor
ring $\R[x]_h$, as has been explained in \cite{S0,S3}. See also \cite{S4}.
   
   Much time is devoted to explaining the properties of Lagrange, Hermite, and other types of
interpolating polynomials in numerical analysis. In teaching this subject, the author has discovered that
many of the formulas and theorems follow directly from the above theorem.  
For rational approximation, we have the following refinement: 
  
   {\bf Rational Approximation Theorem:}  Let $f(x)\in C[a,b]$ and $n$ times differentiable on $(a,b)$.
    Let $u(x)$ and $v(x)$ be polynomials such that $v(0)=1$ and 
   $det(u(x)v(x))< deg(h(x))$, and suppose that $f(x)v(x)-u(x)=0 \ mod(h(x))$. 
   Then
   \[ f(x)=\frac{u(x)}{v(x)}+\frac{1}{n! v(x)}[f(t)v(t)]^{(n)}(c)h(x) \]
for some $c\in (a,b)$.
   
{\bf Proof:} Define the function 
     \[ p(t)=f(t)v(t)-u(t)-[f(x)v(x)-u(x)]\frac{h(t)}{h(x)}  \]
 where $x\in [a,b]$. Clearly, $p(t)=0 \ mod(h(t)(t-x))$.    
Applying the Generalized Rolle's Theorem to $p(t)$, it follows that there exist
   a $c\in (a,b)$ such that 
     \[ f(x)v(x)-u(x)=\frac{1}{n!}\Big(\frac{d}{dt}\Big)^n[f(t)v(t)]_{t\rightarrow c}h(x), \] 
  from which it follows that
     \[ f(x)=\frac{u(x)}{v(x)}+\frac{1}{n! v(x)}[f(t)v(t)]^{(n)}(c)h(x).\]
        
    \hfill $ \framebox{} $    
    
       \section*{Acknowledgements} 
     The author is grateful to Dr. Reyla Navarro, Chairwoman, and
Dr. Guillermo Romero, Vice Rector, of the Universidad
de Las Americas for support for this research. The author is is a member of SNI 14587.

%\begin{address}

{email: garret.sobczyk@udlap.mx

url: http://www.garretstar.com}

%\end{address}

 \end{document}